\documentclass[12pt]{amsart}
\usepackage{graphics}

\newtheorem{theorem}{Theorem}[section]
\newtheorem{lemma}[theorem]{Lemma}

\newtheorem{remarks}[theorem]{Remarks}

\theoremstyle{definition}

\theoremstyle{remark}

\numberwithin{equation}{section}

 \newcommand{\s}{\sigma}
 \newcommand{\D}{\Delta}
 \newcommand{\mcg}{M}
 \newcommand{\jm}{\jmath}
 \def\S{{\Sigma}}
 \def\B{{\rm B}}
 \def\R{{\mathbb{R}}}
 \def\Z{{\mathbb{Z}}}
 \def\CP{{CP}}

\begin{document}

 \title[Noncomplex $4$-manifolds with Lefschetz fibrations]
 {Noncomplex smooth $4$-manifolds with Lefschetz fibrations}
 \author{Mustafa Korkmaz}
 \address{Department of Mathematics, Middle East Technical University,
 06531 Ankara, Turkey} \email{korkmaz@arf.math.metu.edu.tr}
 \subjclass{Primary 57N13, 57N05; Secondary 57R17, 20F38, 20F36}
 \date{\today}
 \keywords{$4$-manifolds, Lefschetz fibrations, Mapping class groups.}

\newenvironment{prooff}{\medskip \par \noindent {\it Proof}\ }{\hfill
$\square$ \medskip \par}
    \def\sqr#1#2{{\vcenter{\hrule height.#2pt
        \hbox{\vrule width.#2pt height#1pt \kern#1pt
            \vrule width.#2pt}\hrule height.#2pt}}}
    \def\square{\mathchoice\sqr67\sqr67\sqr{2.1}6\sqr{1.5}6}
\def\pf#1{\medskip \par \noindent {\it #1.}\ }
\def\endpf{\hfill $\square$ \medskip \par}
\def\demo#1{\medskip \par \noindent {\it #1.}\ }
\def\enddemo{\medskip \par}
\def\qed{~\hfill$\square$}

 \maketitle

 \setcounter{secnumdepth}{1} \setcounter{section}{0}

 \section{Introduction}
 Recently, B. Ozbagci and A. Stipsicz \cite{os} proved that there are infinitely
 many pairwise nonhomeomorphic $4$-manifolds admitting genus-$2$
 Lefschetz fibration over $S^2$ but not carrying any complex structure
 with either orientation. (For the definition of Lefschetz fibration, see
 \cite{gs}.) Their result depends on a relation in the mapping
 class group of a closed orientable surface of genus $2$.
 This relation with eight right Dehn
 twists was discovered by Y. Matsumoto \cite{m} by a computer calculation, and
 it is the global monodromy of a Lefschetz fibration
 $T^2\times S^2\# 4\overline{\CP^2}\to S^2$, where $S^2$ is the $2$-sphere and
 $T^2$ is the $2$-torus.

 In this paper, we generalize Matsumoto's relation to higher genus orientable surfaces.
 We find a relation involving $2g+4$ (resp., $2g+10$) Dehn twists when
 the genus of the surface is even (resp., odd). Following the method of Ozbagci and
 Stipsicz, for every positive integer $n$, we obtain a $4$-manifold $X_n$ admitting
 a genus-$g$ Lefschetz fibration such that the fundamental group of $X_n$
 is isomorphic to $\Z\oplus \Z_n$ for every $g\geq 2$. We then deduce that
 the $4$-manifold $X_n$ does not admit any complex structure.
 This is the main result of this paper.

 \begin{theorem}
 For every $g\geq 2$, there are infinitely many pairwise nonhomeomorphic
 $4$-manifolds that admit genus-$g$ Lefschetz fibrations over $S^2$ but
 do not carry any complex structure with either orientation.
 \label{thm1.1}
 \end{theorem}

 Our relation in the mapping class group given by Theorem
 \ref{relation} also shows that the minimal number of
 singular fibers in a nontrivial genus-$g$ Lefschetz fibration over
 $S^2$ is less than or equal to $2g+4$ (resp., $2g+10$) if $g$ is even (resp., odd).
 This result was also obtained independently by C. Cadavid \cite{c}.
 By definition, a Lefschetz fibration is nontrivial if it admits singular fibers.
 Stipsicz proved in
 \cite{s2} that this minimal number is in fact $2g+4$ (resp., $2g+10$)
 if $g$ is even (resp., odd) and greater than or equal to $6$
 (resp., greater than or equal to $15$) among all $4$-manifolds with
 $b_2^+=1$.  See \cite{ko} and \cite{s}
 for the other results related to this minimal number.

 Here is how we obtain our relation in the mapping class group. Let $\S_g$ be
 a closed connected orientable surface of genus $g$.
 The hyperelliptic
 mapping class group of $\S_g$ is a quotient of the braid group $\B_{2g+2}$ on $2g+2$
 strings. The quotient of the hyperelliptic mapping class group with the
 cyclic subgroup of order $2$ generated by the hyperelliptic involution
 is isomorphic to the mapping class group of a sphere with $2g+2$ holes.
 The hyperelliptic mapping class group is equal to the mapping
 class group when $g=2$. Using these facts, we lift Matsumoto's
 relation to the braid group $\B_6$ and generalize it to a relation in
 $\B_{2g+2}$, although we do not say so explicitly.
 We then project it to the surface $\S_g$ to get our relation
 in the mapping class group of $\S_g$.

 For each positive integer $n$, by considering a product of conjugates of
 our relation with appropriate mapping classes, we obtain
 a relation in the mapping class group of $\S_g$ so that the
 fundamental group of the corresponding symplectic $4$-manifold
 $X_n$ is isomorphic to $\Z\oplus \Z_n$. It follows from
 \cite[proof of Theorem $1.3$]{os} that a symplectic manifold with fundamental
 group $\Z\oplus \Z_n$ admits no complex structures.

 In the last section, we determine the diffeomorphism type of the $4$-manifold
 $X$ admitting a genus-$g$ Lefschetz fibration over $S^2$
 corresponding to our relation given by Theorem \ref{relation}.
 If $g$ is even, then $X$ is diffeomorphic to
 $\S_{g/2}\times S^2\# 4\overline{\CP^2}$. For this, we
 use a result of Stipsicz asserting
 that the only $4$-manifold with $b_2^+=1$ which admits a
 (relatively minimal) Lefschetz fibration with $2g+4$ vanishing cycles
 is $\S_{g/2}\times S^2\# 4\overline{\CP^2}$.

 \section{Braid groups}
 \label{section2}

 The braid group  $\B_{2g+2}$ on $2g+2$ strings admits a
 presentation with generators
 $\s_1,\s_2, \ldots,\s_{2g+1}$ and relations
 $$ \s_i\s_j=\s_j\s_i\,\, \mbox{\rm  if }\,\,|i-j|\geq 2$$
    and
 $$ \s_i\s_{i+1}\s_i=\s_{i+1}\s_{i}\s_{i+1}\,\, \mbox{\rm if }\,\,1\leq i \leq 2g.$$
 The subgroup of $\B_{2g+2}$ generated by
 $\s_1,\s_2, \ldots,\s_h$ is isomorphic to $\B_{h+1}$.
 We identify this subgroup with $\B_{h+1}$.

 In the group $\B_{2g+2}$, let us define the words $\D_k=\s_1\s_2\cdots \s_k$ and
 $\bar{\D}_k=\s_k\cdots \s_2\s_1$ for each $k=1,2,\ldots,2g+1$. We
 define $\D_0=1$ for the convention. For each $k=0,1,\ldots,g,$ let us
 also define
 \[\beta_k=\bar{\D}_k \D_{2g+1-k} \D_{2g-k}^{-1}
         \bar{\D}_k^{-1}
 \]
 and
 \[
 \beta=\bar{\D}_{g}^{g+1}.
 \]
 We have, for example, $\beta_0=\D_{2g+1}\D_{2g}^{-1}$,
 $\beta_1=\bar{\D}_1 \D_{2g} \D_{2g-1}^{-1}\bar{\D}_1^{-1}$, and
 $\beta_g=\bar{\D}_g \D_{g+1} \D_{g}^{-1}\bar{\D}_g^{-1}$. Note that
 $\beta_k$ is the conjugate of $\s_{2g+1-k}$ with the element $\bar{\D}_{k}
 \D_{2g-k}$.

 The following lemma follows easily from the defining relations of the
 braid group.

 \bigskip

 \begin{lemma} The following relations hold in the group $\B_{2g+2}$.
 \begin{enumerate}
 \item[$(a)$] $\s_k\D_m=\D_m \s_{k-1}$ and
    $\s_k^{-1}\D_m =\D_m \s_{k-1}^{-1}$ if $1<k\leq m$;
 \item[$(b)$] $\s_k\bar{\D}_m=\bar{\D}_m \s_{k+1}$ and
    $\s_k^{-1}\bar{\D}_m=\bar{\D}_m \s_{k+1}^{-1}$ if $1\leq k<m$;
 \item[$(c)$] $\s_k\D_m=\D_m \s_{k}$ and
    $\s_k\bar{\D}_m=\bar{\D}_m \s_k$ if $k>m+1$; and
 \item[$(d)$] $\D_g^k=\D_{g-1}\D_g^{k-1} \s_{g-k+1}$ and
    $\bar{\D}_g^k=\s_{g-k+1}\bar{\D}_g^{k-1}\bar{\D}_{g-1}$ if $1\leq k\leq g$.
 \end{enumerate}
 \label{lemma1}
 \end{lemma}

 \bigskip

 \begin{lemma}
 In the braid group $\B_{2g+2}$, we have the following:
     \begin{enumerate}
         \item[$(a)$] The element $\beta$ is equal to
         $\bar{\D}_g \D_g \bar{\D}_{g-1}^{g}$;
    and
           \item[$(b)$] The element $\bar{\D}_g \D_g$ is in the centralizer
    of $\B_g$; in particular, it commutes with $\bar{\D}_{g-1}^{g}$.
 \end{enumerate}
 \label{lemma2}
 \end{lemma}

 \begin{proof}
 We claim that $\bar{\D}_{g}^{g}=\D_k\bar{\D}_g^{g-k}
 \bar{\D}_{g-1}^{k}$ for every $0\leq k\leq g$.
 First of all, the claim holds trivially for $k=0$.
 Suppose by induction that $\bar{\D}_{g}^{g}=\D_k\bar{\D}_g^{g-k}
 \bar{\D}_{g-1}^{k}$. By Lemma \ref{lemma1} $(d)$, we have
 $\bar{\D}_g^{g-k}=\s_{k+1}\bar{\D}_g^{g-k-1} \bar{\D}_{g-1}$. Then
 we have
 \begin{eqnarray*}
 \bar{\D}_{g}^{g} &=& \D_k\bar{\D}_g^{g-k} \bar{\D}_{g-1}^{k}\\
 &=& \D_k \s_{k+1}\bar{\D}_g^{g-k-1} \bar{\D}_{g-1}  \bar{\D}_{g-1}^{k}\\
 &=& \D_{k+1}\bar{\D}_g^{g-(k+1)} \bar{\D}_{g-1}^{k+1}.
 \end{eqnarray*}
 In particular, $\bar{\D}_{g}^{g}=\D_g\bar{\D}_{g-1}^{g}$. The proof
 of $(a)$ follows.

 For $k< g$, by Lemma \ref{lemma1}, we have \( \bar{\D}_g \D_g \s_k
 =\bar{\D}_g\s_{k+1} \D_g = \s_k\bar{\D}_g \D_g. \) Now the proof of $(b)$
 follows.
 \end{proof}

 \bigskip
 \begin{lemma}
 For every $0\leq k\leq g-1$, we have
 \[
 \D_{2g-k}^{-1}   \bar{\D}_k^{-1} \bar{\D}_{k+1}
 \D_{2g-k}  \bar{\D}_{k+1}  \D_{2g-k-1} =\bar{\D}_k
 \D_{2g-k} \gamma_k,
 \]
 where $\gamma_k=\bar{\D}_k \D_{2g-k-1} \D_{2g-k-2}^{-1} \bar{\D}_k^{-1}$.
 \label{lemma3}
 \end{lemma}

 \begin{proof}
 We use Lemma \ref{lemma1} several times:
 \begin{eqnarray*}
 \D_{2g-k}^{-1}   & \!\!\!\! \bar{\D}_k^{-1}  &\!\!\!\!
    \bar{\D}_{k+1}
    \D_{2g-k}  \bar{\D}_{k+1}  \D_{2g-k-1}  \\
 &&=\s^{-1}_{2g-k} \cdots \s^{-1}_{k+2} \s^{-1}_{k+1} \D_{k}^{-1}
    \bar{\D}_{k}^{-1} (\bar{\D}_{k+1} {\D}_{2g-k} \bar{\D}_{k+1}
    {\D}_{2g-k-1}) \\
 &&=\s^{-1}_{2g-k}\cdots\s^{-1}_{k+2}
    \s^{-1}_{k+1} (\bar{\D}_{k+1} {\D}_{2g-k} \bar{\D}_{k+1}
    {\D}_{2g-k-1}) \D_{k}^{-1} \bar{\D}_{k}^{-1}  \\
 &&=\s^{-1}_{2g-k}\cdots \s^{-1}_{k+3}
    \s^{-1}_{k+2} \bar{\D}_{k} {\D}_{2g-k} \bar{\D}_{k+1}
    {\D}_{2g-k-1} \D_{k}^{-1} \bar{\D}_{k}^{-1}  \\
 &&=\bar{\D}_{k}\s^{-1}_{2g-k} \cdots \s^{-1}_{k+3} \s^{-1}_{k+2}
    {\D}_{2g-k}\bar{\D}_{k+1} {\D}_{2g-k-1} \D_{k}^{-1}
    \bar{\D}_{k}^{-1}  \\
 &&=\bar{\D}_{k} {\D}_{2g-k} \s^{-1}_{2g-k-1} \cdots \s^{-1}_{k+2}
    \s^{-1}_{k+1} \bar{\D}_{k+1} {\D}_{2g-k-1} \D_{k}^{-1}
    \bar{\D}_{k}^{-1} \\
 &&=\bar{\D}_{k} {\D}_{2g-k}
    \s^{-1}_{2g-k-1} \cdots \s^{-1}_{k+3} \s^{-1}_{k+2} \bar{\D}_{k}
    {\D}_{2g-k-1} \D_{k}^{-1} \bar{\D}_{k}^{-1} \\
 &&=\bar{\D}_{k}{\D}_{2g-k}\bar{\D}_{k} {\D}_{2g-k-1}
    \s^{-1}_{2g-k-2} \cdots \s^{-1}_{k+2} \s^{-1}_{k+1} \D_{k}^{-1}
    \bar{\D}_{k}^{-1} \\
 &&=\bar{\D}_{k} {\D}_{2g-k}\bar{\D}_{k}
    {\D}_{2g-k-1} \D_{2g-k-2}^{-1} \bar{\D}_{k}^{-1} \\
 &&=\bar{\D}_{k}{\D}_{2g-k}\gamma_k.
 \end{eqnarray*}
 \end{proof}

 The main result of this section is the next theorem.

 \begin{theorem}
 \label{thm3}
 In the braid group $\B_{2g+2}$, we have the relation
 \[
 \beta_0\beta_1\beta_2\cdots \beta_g\beta^2=\D_{2g+1}\D_{2g}
                                            \cdots\D_3\D_2\D_1.\]
 \end{theorem}

 \begin{proof}
 Recall that, for any $h\leq 2g+2$, we identify the group $\B_{h}$ with the
 subgroup of $\B_{2g+2}$ generated by the elements $\s_1,\s_2,\ldots,\s_{h-1}$.

 The proof of the theorem is by induction on $g$.
 Suppose that $g=0$. In the group $\B_2$,
 $\beta_0=\D_1\D_0^{-1}$ and $\beta=\D_0$. Thus
 $\beta_0\beta^2=\D_1$. Hence, the conclusion of the theorem holds for $g=0$.

 In the subgroup $\B_{2g}$ of $\B_{2g+2}$, let us define
 \[ \gamma_k=\bar{\D}_k\D_{2g-1-k}\D_{2g-2-k}^{-1}\bar{\D}_k^{-1},\,
            0\leq k\leq g-1 \]
 and
 \[ \gamma=\bar{\D}_{g-1}^g.\]
 Then, by the induction hypothesis
 \[ \gamma_0\gamma_1\gamma_2\cdots\gamma_{g-1}\gamma^2
 =\D_{2g-1}\D_{2g-2}\cdots \D_{3}\D_{2}\D_{1}.
 \]
 Let us also define $\gamma_g=1$ for the convention.

 In the group $\B_{2g+2}$, we claim that
 \[
 \beta_k\beta_{k+1}\cdots\beta_g\beta^2=
    \bar{\D}_{k} \D_{2g+1-k}\bar{\D}_{k}\D_{2g-k}
 \gamma_k\gamma_{k+1}\cdots\gamma_{g-1}\gamma_{g}\gamma^2.
 \]
 The proof of this claim is by induction on $g-k$. We start with the
 following computation:
 \begin{eqnarray*} \beta_g\beta^2
        &=& \bar{\D}_g \D_{g+1}\D_{g}^{-1} \bar{\D}_{g}^{-1}
        (\bar{\D}_g \D_g \bar{\D}_{g-1}^{g})^2\\
        &=& \bar{\D}_g \D_{g+1} \bar{\D}_{g-1}^{g} \bar{\D}_g \D_g \bar{\D}_{g-1}^{g}\\
        &=& \bar{\D}_g \D_{g+1} \bar{\D}_g \D_g
        (\bar{\D}_{g-1}^{g})^2\\
        &=& \bar{\D}_g \D_{g+1} \bar{\D}_g \D_g
        \gamma_g\gamma^2.
 \end{eqnarray*}
 Hence, the claim holds for $k=g$. Suppose inductively that
 \[
 \beta_{k+1}\beta_{k+2}\cdots\beta_g\beta^2=
        \bar{\D}_{k+1}\D_{2g-k}\bar{\D}_{k+1}\D_{2g-k-1}
        \gamma_{k+1}\gamma_{k+2}\cdots\gamma_{g}\gamma^2.
 \]
 Then by Lemma \ref{lemma3} we get
 \begin{eqnarray*}
 &&\hspace*{-1cm}\beta_k\beta_{k+1}\cdots \beta_g\beta^2\\
  &&=\bar{\D}_{k} \D_{2g+1-k} (\D_{2g-k}^{-1}\bar{\D}_{k}^{-1}
     \bar{\D}_{k+1}\D_{2g-k}\bar{\D}_{k+1}\D_{2g-k-1})
  \gamma_{k+1}\gamma_{k+2}\cdots\gamma_{g}\gamma^2\\
  &&=\bar{\D}_{k} \D_{2g+1-k}\bar{\D}_{k} \D_{2g-k}
     \gamma_{k}\gamma_{k+1}\cdots\gamma_{g}\gamma^2.
 \end{eqnarray*}
 Hence, the claim is proved. For $k=0$, in particular, we obtain
 \begin{eqnarray*}
 \beta_0\beta_1\beta_2\cdots\beta_g\beta^2
    &=&\bar{\D}_{0} \D_{2g+1}\bar{\D}_{0} \D_{2g}
    \gamma_{0}\gamma_{1}\cdots\gamma_{g}\gamma^2\\
    &=&\D_{2g+1} \D_{2g}
    \gamma_{0}\gamma_{1}\cdots\gamma_{g}\gamma^2\\
    &=&\D_{2g+1} \D_{2g}\D_{2g-1}\cdots\D_3\D_2\D_1.
 \end{eqnarray*}

 This finishes the proof of the theorem.
 \end{proof}

 \section{Mapping class groups}
 \label{section3}

 Let $\S_g$ be a closed connected orientable surface of genus $g$ embedded in $\R^3$
 such that it is invariant under the involution
 $J(x,y,z)=(-x,y,-z)$ (cf. Fig. \ref{sekil5}).
 Notice that $J$ is the rotation about $y$-axis by $\pi$.
 We orient $\S_g$ so that the unit
 normal vectors are pointing outward. Let $\jm$ be the isotopy class
 of $J$. Let us denote
 by $\mcg_g$ the mapping class group of $\S_g$ which is the group of isotopy
 classes of orientation-preserving diffeomorphisms of $\S_g$. The
 hyperelliptic mapping class group is defined as the centralizer
 $C_{\mcg_g}(\jm)$ of $\jm$ in $\mcg_g$, the subgroup consisting of those
 the mapping classes that commute with $\jm$.

 Throughout this paper, we use functional notation. That is, for
 any two mapping classes $f$ and $g$, the multiplication $fg$ means that
 $g$ is applied first.

 For a simple closed curve $a$ on the oriented surface $\S_g$,
 by the abuse of notation, a right Dehn twist about $a$ and its isotopy class
 is denoted by $t_a$.

 Let us consider the simple closed curves $A_1,A_2,\ldots,A_{2g+1}$ on
 $\S_g$ defined as follows: $A_{2k}=b_k$, $A_1=a_1$,
 $A_{2k-1}=a_ka_{k+1}^{-1}$, and $A_{2g+1}=a_g$, where $a_k$ and
 $b_k$ are the curves shown  in Fig. \ref{sekil5}. Let $t_k$ denote the right Dehn
 twist about $A_k$. The simple closed curves $A_k$ are invariant under $J$. It follows
 that Dehn twists $t_1,t_2,\dots,t_{2g+1}$ commute with $\jm$. Hence, they
 are contained in the hyperelliptic mapping class group.

 The involution $J$ has $2g+2$ fixed points. Hence, we have a branched
 covering $p:\S_g\to S^2$ branching over $2g+2$ points. Let us
 mark these $2g+2$ points on $S^2$, and let $\S_{0,2g+2}$ be the resulting
 surface.  Notice that the interior of each $p(A_k)$ is an embedded
 arc on $\S_{0,2g+2}$ connecting two distinct marked points
 and that it is disjoint from $p(A_l)$ for $k\neq l$. Let us denote by $w_k$ the isotopy
 class of a right half twist about $p(A_k)$. Thus, if we orient the arc
 $p(A_k)$ arbitrarily, then $w_k(p(A_k))$ (defined up to isotopy) is isotopic to the arc
 $p(A_k)^{-1}$. Therefore, $w_k^2$ is the right Dehn twist
 about the boundary component of a regular neighborhood of $p(A_k)$. It is
 well-known that the half twists $w_1,w_2,\ldots,w_{2g+1}$ generate the mapping class
 group ${\mcg}_{0,2g+2}$ of $\S_{0,2g+2}$ (cf. \cite[Theorem
 $4.5$]{b}). Here, the group ${\mcg}_{0,2g+2}$ is defined to be the
 group of the isotopy classes of the orientation-preserving
 diffeomorphisms of $\S_{0,2g+2}$ that
 preserve the marked points setwise. The isotopies are assumed
 to fix each marked point.

 \begin{theorem}
 The hyperelliptic mapping class group $C_{\mcg_g}(\jm)$ is generated by
 the Dehn twists $t_1,t_2,\ldots,t_{2g+1}$, the function given by $\Psi (t_k)=w_k$
 on the generators defines a surjective homomorphism
 \[
 \Psi :\, C_{\mcg_g}(\jm)\to {\mcg}_{0,2g+2},
 \]
 and the kernel of $\Psi$ is $\langle \jm \rangle$, which is a subgroup of order $2$.
 \label{hmcg}
 \end{theorem}

 Theorem \ref{hmcg} was proved by J. Birman and H. Hilden \cite{bh}. They
 also obtained a presentation of the hyperelliptic mapping class group.
 Since $t_it_j=t_jt_i$ for $|i-j|\geq 2$ and $t_{i}t_{i+1}t_{i}=t_{i+1}t_{i}t_{i+1}$
 in the group $C_{\mcg_g}(\jm)$, the fact that $t_1,t_2,\ldots,t_{2g+1}$
 generate $C_{\mcg_g}(\jm)$ implies that
 $\s_k \mapsto t_k$ defines a surjective homomorphism
 $B_{2g+2}\to C_{\mcg_g}(\jm)$.

\begin{figure}[hbt]
 \begin{center}
    \includegraphics{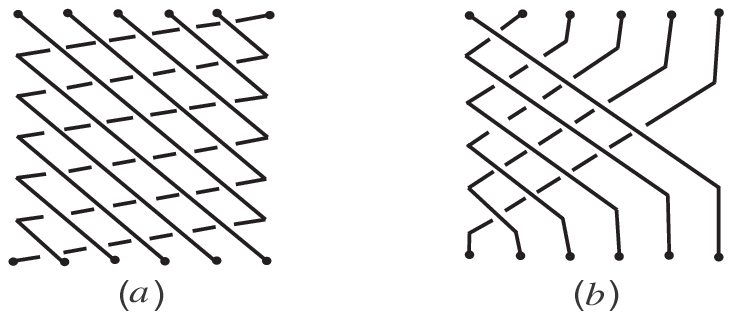}
  \caption{The words $(\s_5\s_4\s_3\s_2\s_1)^6$ and $\D$ in $\B_6$.}
  \label{sekil2}
   \end{center}
 \end{figure}

 The following lemma is easy to prove (cf. Fig. \ref{sekil2} $(a)$).

 \begin{lemma}
 In the group ${\mcg}_{0,2g+2}$, the element
 $(w_k\cdots w_2w_1)^{k+1}$ is equal to the right Dehn twist about the
 boundary component of a regular neighborhood of $p(A_1)\cup
 p(A_2)\cup\cdots\cup p(A_k)$.
 \label{delta2=1}
 \end{lemma}

\begin{figure}[hbt]
 \begin{center}
    \includegraphics{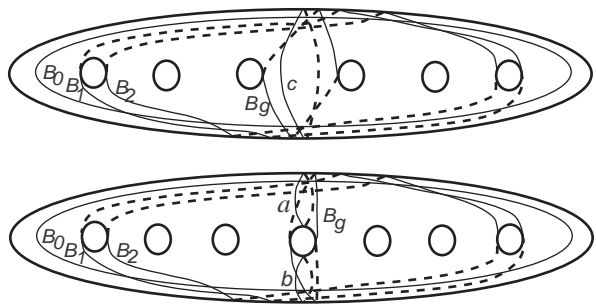}
  \caption{The simple closed curves $B_k$, $a$, $b$ and $c$.}
  \label{sekil3}
   \end{center}
 \end{figure}

 In order to state the main result of this section, let us consider the
 simple closed curves $B_k,a,b,$ and $c$ illustrated in Fig. \ref{sekil3}.
 Note that $a$ and $b$ are defined for odd $g$, and $c$ is defined for
 even $g$.

 \begin{lemma}
 The following relations hold in the mapping class group:
 \begin{enumerate}
  \item[$(a)$] $t_c=(t_g \cdots t_2 t_1)^{2(g+1)}$ if $g$ is even; and
  \item[$(b)$] $t_a t_b= (t_g \cdots t_2 t_1)^{g+1}$ if $g$ is odd.
 \end{enumerate}
 \label{lemma2.4}
 \end{lemma}

\begin{figure}[hbt]
 \begin{center}
    \includegraphics{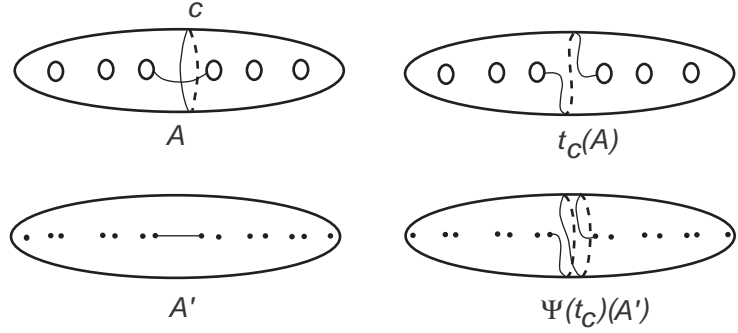}
  \caption{Projection of the Dehn twist $t_c$.}
  \label{sekil4}
   \end{center}
 \end{figure}

 \begin{proof}
 Suppose that $g$ is even. Consider the branched covering
 $p:\S_g\to S^2$.
 Notice that $p(c)$ is a simple closed curve on $\S_{0,2g+2}$.
 The projection $\Psi (t_c)$ of $t_c$ to $\S_{0,2g+2}$ is the square of the Dehn twist
 about $p(c)$. This can be seen geometrically as follows.
 Consider the arcs $p(A_1),p(A_2),\ldots ,p(A_{2g+1})$. Since the surface
 obtained by cutting $\S_{0,2g+2}$ along these arcs is a disc without any marked points
 in the interior, in order to show that $\Psi (t_c)=t_{p(c)}^2$,
 it is enough to check that the actions of $\Psi
 (t_c)$ and $t_{p(c)}^2$ on these arcs are the same (up to isotopy). To see the action of
 $\Psi (t_c)$ on an arc $A'$, lift $A'$ to $\S_g$, apply $t_c$,
 and then project  it down to $\S_{0,2g+2}$ (cf. Fig. \ref{sekil4}).
 Since $t_{p(c)}=(w_g \cdots w_2w_1)^{g+1}$ by Lemma \ref{delta2=1}, we conclude that
 $\Psi (t_c)=\Psi ((t_g\cdots t_2t_1)^{2(g+1)})$. Hence,
 $(t_g\cdots t_2
 t_1)^{2(g+1)}$ is equal to either $t_c$ or $\jm t_c$. We rule out the latter
 possibility as follows. It is easily
 checked that $(t_g\cdots t_2t_1)^{2(g+1)}$ acts trivially on the first
 homology group $H_1(\S_g;\Z)$. The action of $t_c$ is also trivial
 since $c$ is null homologous. On the other hand, $\jm$ acts as the minus
 identity on $H_1(\S_g;\Z)$.

 The part $(b)$ is proved similarly.
 \end{proof}

 \begin{theorem}
 In the mapping class group $\mcg_g$, the following
 relations between right Dehn twists hold:
 \begin{itemize}
 \item[$(a)$] $(t_{B_0}t_{B_1}t_{B_2}\cdots t_{B_g}t_c)^2=1$ if $g$ is even;
 \item[$(b)$] $(t_{B_0}t_{B_1}t_{B_2}\cdots t_{B_g}t_a^2t_b^2)^2=1$  if $g$ is
 odd.
 \end{itemize}
 \label{relation}
 \end{theorem}

 \begin{proof}
 In the mapping class group $\mcg_g$ of $\S_g$, for each
 $k=0,1,2,\ldots, g$, we define $\D_k,\bar{\D}_k,\beta_k,$ and $\beta$ as in
 Section \ref{section2} by replacing $\s_i$ by $t_i$. Recall that $t_i$ is
 the (right) Dehn twist about the simple closed curve $A_i$. Hence,
 \[
 \beta_k=(\bar{\D}_k\D_{2g-k}) t_{2g+1-k} (\bar{\D}_k\D_{2g-k})^{-1}.
 \]
 It is easy to see that $\bar{\D}_i\D_{2g-i}(A_{2g+1-i})=B_i$.
 Since $ft_ef^{-1}=t_{f(e)}$ for any $f\in \mcg_g$ and for any
 simple closed curve $e$, we conclude that
 $\beta_k=t_{B_k}$. Also, by Lemma \ref{lemma2.4}, $\beta^2=t_c$ if $g$ is even
 and $\beta^2=t_a^2t_b^2$ if $g$ is odd. Let us define the word
 \[
 W=\left\{
 \begin{array}{ll}
 (t_{B_0}t_{B_1}t_{B_2}\cdots t_{B_g}t_c)^2 & {\rm if}\,\,g\,\,
 {\rm is\,\,even,}\\
 (t_{B_0}t_{B_1}t_{B_2}\cdots t_{B_g}t_a^2t_b^2)^2 & {\rm if}\,\,g\,\,
 {\rm is\,\, odd}.
 \end{array}
 \right.
 \]
 Hence, $W=(\beta_0\beta_1\cdots\beta_g\beta^2)^2$. Let
 $\D=\D_{2g+1}\D_{2g} \cdots \D_2\D_1$. Since $t_it_j=t_jt_i$ for $|i-j|\geq 2$
 and $t_it_{i+1}t_i=t_{i+1}t_it_{i+1}$ by Theorem \ref{thm3}, we obtain $W=\D^2$.
 Now the element $\Psi (\D)$ is of order $2$ in $\mcg_{0,2g+2}$ (cf. Fig. \ref{sekil2}
 $(b)$ for $g=6$), where $\Psi$ is the epimorphism in Theorem \ref{hmcg}.
 Hence, either $\D^2=1$ or $\D^2=\jm$. It is easy to verify
 that $\D^2$ acts trivially on the first homology group $H_1(\S_g;{\bf
 Z})$ but $\jm$ does not. Therefore, $W=1$.

 This finishes the proof of the theorem.
 \end{proof}

 \section{Noncomplex genus-$g$ Lefschetz fibrations}
 \label{section4}

\begin{figure}[hbt]
 \begin{center}
    \includegraphics{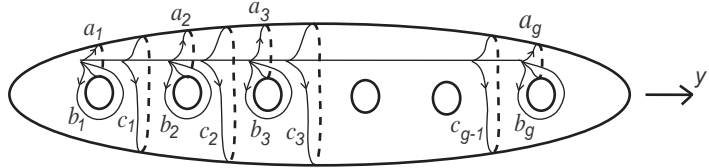}
  \caption{Generators of the fundamental group.}
  \label{sekil5}
   \end{center}
 \end{figure}

 In this section, we prove the main result of this paper. We assume from
 now on that  $g\geq 2$.
 We construct a $4$-manifold $X_n$ a admitting genus-$g$ Lefschetz
 fibration with fundamental group $\Z\oplus \Z_n$ for every
 positive integer $n$. Then we conclude the main result of this paper
  from the proof of the
 main result of \cite{os}. As the model for a closed connected oriented
 surface $\S_g$, we will consider the one embedded in $\R^3$
 as shown in Fig. \ref{sekil5}.

 For any two elements $x$ and $y$ in a group, we denote $yxy^{-1}$ and
 $xyx^{-1}y^{-1}$ by $x^y$ and $[x,y]$, respectively.

 Let us consider the word $W$ in $\mcg_g$ defined by
 \[
 W=\left\{
 \begin{array}{ll}
 (t_{B_0}t_{B_1}t_{B_2}\cdots t_{B_g}t_c)^2 & {\rm if}\,\,g\,\,
 {\rm is\,\,even,}\\
 (t_{B_0}t_{B_1}t_{B_2}\cdots t_{B_g}t_a^2t_b^2)^2 & {\rm if}\,\,g\,\,
 {\rm is\,\, odd}.
 \end{array}
 \right.
 \]
 By Theorem \ref{relation}, $W=1$ in $\mcg_g$. Since the conjugation of a
 right Dehn twist with an element of $\mcg_g$ is again a right Dehn twist,
 the word $W^{f}$ is a product of right Dehn twists for any mapping class $f$.
 For each positive integer $n$, we define
 \[
 W_n=\left\{
 \begin{array}{ll}
 WW^{t_{a_1}} W^{t_{a_2}} \cdots W^{t_{a_{r-1}}} W^{t_{a_r}^n} W^{t_{b_{r+2}}}
 W^{t_{b_{r+3}}} \cdots W^{t_{b_g}} & {\rm if}\,\,g=2r,\\
 WW^{t_3^n} W^{t_5} W^{t_7}  \cdots W^{t_{2r+1}}  W^{t_{b_{r+2}}}
 W^{t_{b_{r+3}}} \cdots W^{t_{b_g}} & {\rm if}\,\,g=2r+1,
 \end{array}
 \right.
 \]
 where $a_i$ and $b_i$ are simple closed curves given in Fig. \ref{sekil5}
 (considered up to isotopy). Note that $W_n$ is a product of $g(2g+4)$
 (resp., $g(2g+10)$) right Dehn twists
 if $g$ is even (odd).

 The word $W_n$ is equal to the identity in the mapping class group $\mcg_g$.
 Let $X$ and $X_n$ be the smooth $4$-manifolds that admit the genus-$g$
 Lefschetz fibrations over $S^2$ whose global monodromies are $W$ and
 $W_n$, respectively. Thus, $X_n$ is  the fiber sum of
 $g$ copies of $X$.

 \begin{theorem}
 The fundamental group $\pi_1(X_n)$ of $X_n$ is isomorphic to $\Z\oplus\Z_n$.
 \label{pi1}
 \end{theorem}

 \begin{proof}
 Let $a_k$ and $b_k$ be the standard generators of
 $\pi_1(\S_g)$ illustrated in Fig. \ref{sekil5}.  By the theory of
 Lefschetz fibrations, $\pi_1(X_n)$ is isomorphic to the quotient
 of $\pi_1(\S_g)$ by the normal subgroup generated by the vanishing cycles.

 Suppose first that $g=2r$. It is easy to check that up to conjugation
 the following equalities hold in $\pi_1(\S_g)$:
 \begin{itemize}
  \item $B_0=b_1b_2\cdots b_g;$
  \item $B_{2k-1}=a_{k}b_kb_{k+1}\cdots b_{g+1-k}
                     c_{g+1-k}a_{g+1-k}$, $1\leq k\leq r$;
  \item $B_{2k}=a_{k}b_{k+1}b_{k+2}\cdots b_{g-k}
                     c_{g-k}a_{g+1-k}$, $1\leq k\leq r-1$;
  \item $B_{g}=B_{2r}=a_{r}c_{r}a_{r+1}$;
  \item $c=c_{r}=[a_{1},b_{1}][a_{2},b_{2}]\cdots [a_{r},b_{r}]$.
 \end{itemize}

 The vanishing cycles corresponding to $W^{t_{a_k}}$ are the set
 \[
 \{ a_kB_0,\ldots,a_kB_{2k-1},B_{2k},\ldots,B_g,c\},\,\,1\leq k\leq r-1.\]
 Similarly, the vanishing cycles corresponding to $W^{t_{a_r}^n}$ and
 $W^{t_{b_{g+1-k}}}$ are
 \[
 \{ a_r^nB_0,\ldots,a_r^nB_{g-1},B_g,c\}\] and
 \[
 \{ B_0,\ldots,B_{2k-2},b_{g+1-k}^{-1}B_{2k-1},b_{g+1-k}^{-1}B_{2k},
 B_{2k+1},\ldots,B_g,c\},\,\,1\leq k\leq r-1,\]
 respectively. It follows that the
 fundamental group of $X_n$ has a presentation with generators
 $a_1,b_1,a_2,b_2,\ldots,a_g,b_g$ and relations
 \begin{itemize}
  \item $\Pi_{k=1}^g [a_k,b_k]=1$;
  \item $B_0=B_1=B_2=\cdots =B_g=c=1$;
  \item $a_1=a_2=\cdots
  =a_{r-1}=a_r^n=b_{r+2}=b_{r+3}=\cdots=b_{g}=1$.
 \end{itemize}
 It is easy to see that this presentation is equivalent to the
 presentation with generators
 $a_r,b_r$ and relations $a_r^n=[a_r,b_r]=1$.

 Suppose now that $g=2r+1$. A similar argument as in the case of
 even $g$ shows that the fundamental group of $X_n$ has a presentation
 with generators
 $a_1,b_1,a_2,b_2,\ldots,a_g,b_g$ and relations
 \begin{itemize}
  \item $\Pi_{k=1}^g [a_k,b_k]=1$;
  \item $B_0=B_1=B_2=\cdots =B_g=a=b=1$;
  \item $(a_2a_1^{-1})^n=a_3a_2^{-1}=a_4a_3^{-1}=\cdots
  =a_{r+1}a_r^{-1}=b_{r+2}=b_{r+3}=\cdots=b_{g}=1$.
 \end{itemize}
 Since $a=a_{r+1}$, this presentation is equivalent to the presentation
 with generators $a_1,b_1$ and relations $a_1^n=[a_1,b_1]=1$.

 This completes the proof of the theorem.
 \end{proof}

 The following theorem can be concluded from \cite[proof of
 Theorem $1.3$]{os}.

 \begin{theorem} \label{thm4.2}
 Let $M$ be an orientable $4$-manifold such that $b_2^+(M) \geq 1$ and
 $\pi_1(M)=\Z\oplus \Z_n$. Then $M$ does not carry any complex structure.
 \end{theorem}

 \begin{prooff}{\em of Theorem~\ref{thm1.1}.}
 The smooth $4$-manifold $X_n$ admits a genus-$g$ Lefschetz fibration for
 every positive integer $n$. Since $X_n$ is symplectic by a result of
 R. Gompf \cite{gs}, we have $b_2^+(X_n)\geq 1$. We showed above that
 $\pi_1(X_n)=\Z\oplus \Z_n$. Hence, the manifold $X_n$ is not homeomorphic
 to $X_m$ for $n\neq m$. By Theorem~\ref{thm4.2}, the manifold $X_n$
 and the manifold obtained from it by reversing the orientation do not admit
 any complex structure.
 \end{prooff}

 \section{The $4$-manifold admitting a Lefschetz fibration with
 global monodromy $W$.}
 \label{section5}

 In this section, we determine the $4$-manifold corresponding
 to the word $W$ given in Section \ref{section3}.

 \subsection{The case of even $g$.}
 Notice that simple closed curves $B_0,B_1,\ldots,B_g$, and $c$
 are invariant under the involution $J$. Hence, the
 genus-$g$ Lefschetz fibration $X\to S^2$ with global monodromy $W$ is
 is hyperelliptic.

 \begin{theorem}[{\cite{m1},\cite{m},\cite{e},\cite{o}}]
 Let $M$ be a $4$-manifold that admits a hyperelliptic Lefschetz
 fibration of genus $g$ over $S^2$. Let $m$ and
 $s=\sum_{h=1}^{[g/2]} s_h$ be the numbers of nonseparating
 and separating vanishing cycles in the global monodromy of this
 fibration, respectively, where $s_h$ denotes the number of
 separating vanishing cycles that separate the genus-$g$ surface
 into two surfaces one of which has genus $h$.  Then the signature of $M$
 is
 \[
 \s(M)=-\frac{g+1}{2g+1}\, m+\sum_{h=1}^{[g/2]} \left(
 \frac{4h(g-h)}{2g+1}-1\right) s_h.
 \]
 \end{theorem}

 \bigskip
 Since there are $2g+4$ vanishing cycles, Euler
 characteristic of $X$ is $\chi(X)=2(2-2g)+2g+4=8-2g$. There are only two
 separating vanishing cycles, and they bound a surface of genus $g/2$
 on both sides. Hence, the signature of $X$ is
  \[
  \s(X)=-\frac{g+1}{2g+1}\, (2g+2) + \left(
   \frac{4\frac{g}{2}(g-\frac{g}{2})}{2g+1}-1\right) 2=-4.
 \]
 The group $\pi_1(X)$ has a presentation with generators
 $a_1,b_1,\ldots,a_g,b_g$ and relations
 $\Pi_{k=1}^g [a_k,b_k]=B_0=B_1=B_2=\cdots =B_g=c=1$. It is now easy to see that
 $H_1(X;\Z)=\Z^g$. In particular, $b_1(X)=g$. It follows from
 $\chi(X)=8-2g$ and $\s(X)=-4$ that  $b_2^+(X)=1$.
 By \cite[Remark $4.5$(a)]{s2}, $X$ is diffeomorphic to
 $\S_{g/2} \times S^2\# 4\overline{\CP^2}$ if $g\geq 6$.

 \subsection{The case of odd $g$.}
 Suppose that $g$ is at least $3$ and odd.
 Let $X\to S^2$ be the genus-$g$ Lefschetz fibration with global
 monodromy $W$. Since there are $2g+10$ singular fibers, the Euler characteristic
 of $X$ is $\chi(X)=2(2-2g)+2g+10=14-2g$. The fundamental group
 $\pi_1(X)$ of $X$ has a presentation with generators
 $a_1,b_1,\ldots,a_g,b_g$ and relations
 $\Pi_{k=1}^g [a_k,b_k]=B_0=B_1=B_2=\cdots =B_g=a=b=1$. It is now easy to see that
 $H_1(X;\Z)=\Z^{g-1}$. In particular, $b_1(X)=g-1$. Hence,
 \[
 14-2g=2-2b_1(X)+b_2(X)=2-2(g-1)+b_2(X);\]
 that is, $b_2(X)=10$.

 The manifold $X$ is symplectic.  Since $1-b_1+b_2^+$
 is even for any symplectic manifold, we conclude that $b_2^+(X)$
 is odd and is between $1$ and $9$. Hence, $b_2^-(X)$
 is also odd and is between $1$ and $9$.

 We now determine the signature of $X$.
 A handlebody decomposition for $X$ is obtained as follows.
 Start with $\S_g\times D^2$, where $D^2$ is the $2$-disc. Its boundary is
 $\S_g\times S^1$. Attach a $2$-handle along each vanishing cycle (by
 counting its multiplicity) with the $-1$ framing relative to the product
 framing. The cores of the first two $2$-handles attached along $a$
 gives us a $(-2)$-sphere $S_1$. Denote the class of $S_1$ in
 $H_2(X;\mathbb{R})$ by $[S_1]$. Similarly, the cores of the second and the third
 $2$-handles, and the third and the fourth $2$-handles give two $(-2)$-spheres
 $S_2$ and $S_3$. Note that $[S_1][S_3]=0$. Orient each $S_i$ so that
 $[S_1][S_2]=[S_2][S_3]=-1$.
 Similarly, the four $2$-handles glued along $b$ give three
 $(-2)$-spheres
 $S_4,S_5,S_6$ with $[S_4][S_5]=[S_5][S_6]=-1$ and $[S_4][S_6]=0$.
 Since $a$ and $b$ are disjoint, $[S_i][S_j]=0$ for $1\leq i\leq 3$ and
 $4\leq j\leq 6$.
 The first handles attached along $a$ and $b$ give a surface $S_7$
 of genus $(g-1)/2$ such that $[S_7]^2=-2$, $[S_7][S_1]=[S_7][S_4]=-1$,
 and $[S_7][S_i]=0$ for $i=2,3,5,6$.

 The homology classes $[S_1],\ldots,[S_7]$ are linearly
 independent. Hence, they form a basis for a subspace $V$ of
 $H_2(X;\mathbb{R})$ of dimension $7$.
 The matrix of the intersection form restricted to $V$ in
 the above basis is the matrix $-A$, where
 \[A=
   \begin{pmatrix}
     2 & 1 & 0 & 0 & 0 & 0 & 1  \\
     1 & 2 & 1 & 0 & 0 & 0 & 0  \\
     0 & 1 & 2 & 0 & 0 & 0 & 0  \\
     0 & 0 & 0 & 2 & 1 & 0 & 1  \\
     0 & 0 & 0 & 1 & 2 & 1 & 0  \\
     0 & 0 & 0 & 0 & 1 & 2 & 0  \\
     1 & 0 & 0 & 1 & 0 & 0 & 2  \\
  \end{pmatrix}
 \]
 It is easily check that the matrix $A$ is positive definite.
 Hence, the restriction of the intersection form to $V$ is
 negative definite.

 On the other hand, we have $[F]\neq 0$, $[F]^2=0$ and $[F]\in
 V^{\perp}$, where $F$ is a generic fiber and
 $V^{\perp}$ is the orthogonal complement
 of $V$. Since the restriction of the intersection form to
 $V^{\perp}$ is nondegenerate, there is a class $[S_8]\in V^{\perp}$
 with $[S_8]^2<0$. Thus, the restriction of the intersection
 form to the $8$-dimensional subspace generated by $[S_1],\ldots,[S_8]$
 is negative definite. Therefore, $b_2^-(X)$ is at least
 $8$, and hence $b_2^-(X)=9$, since it is also odd. Consequently,
 we have $b_2^+(X)=1$ and $\sigma(X)=-8$.

 \bigskip
 \begin{remarks}
 {\rm
 (1) Stipsicz~\cite{s2} points out that when $g$ is odd, the manifold $X$
 that admits a Lefschetz fibration over $S^2$ with global
 monodromy $W$ is diffeomorphic to
 $\S_{(g-1)/2}\times S^2\# 8\overline{\CP^2}$.

 (2) R. Fintushel and R. Stern \cite{fs} constructed infinite classes of simply
 connected homeomorphic but nondiffeomorphic symplectic manifolds, all
 of which admit Lefschetz fibrations of a fixed fiber genus.
 }
 \end{remarks}

 {\bf {Acknowledgments.}}
 I thank Andr\'{a}s I.
 Stipsicz and Sergey Finashin for helpful conversations and for
 suggestions in the computation of the signature in the case
 of odd $g$. I also thank Yildiray Ozan for answering numerous
 questions.

 \bigskip


\begin{thebibliography}{ZZZZZ}

 \bibitem{b}
 J. S. Birman,
 {\em Braids, links and mapping class groups},
 Ann. of Math. Stud. {\bf82}, Princeton Univ. Press, Princeton, 1974.

 \bibitem{bh}
 J. S. Birman and H. M. Hilden,
 "On the mapping class groups of closed surfaces as covering spaces",
 in {\it Advances in the Theory of Riemann surfaces (Stoony Brook, N.Y., 1969)},
 Ann. of Math. Stud. {\bf 66}, Princeton
 Univ. Press, Princeton 1971, 81-115.

 \bibitem{c}
 C. Cadavid,
 {\em A remarkable set of words in the mapping class group},
 Ph.D. dissertation, Univ. of Texas at Austin, 1998.

 \bibitem{e}
 H. Endo,
 {\em Meyer's signature cocycle and hyperelliptic fibrations},
 Math. Ann. {\bf 316} (2000), 237-257.

 \bibitem{fs}
 R. Fintushel and R. Stern,
 {\em Symplectic surfaces in a fixed homology class},
 J. Differential Geom. {\bf 52} (1999), 203-222.

 \bibitem{gs}
 R. E. Gompf and A. I. Stipsicz,
 {\em $4$-manifolds and Kirby calculus},
 Grad. Stud. Math. {\bf 20}, Amer. Math. Soc., Providence, 1999.

 \bibitem{ko}
 M. Korkmaz and B. Ozbagci,
 {\em Minimal number of singular fibers in a Lefschetz fibration},
 to appear in Proc. Amer. Math. Soc.

 \bibitem{m1}
 Y. Matsumoto,
 {\em On $4$-manifolds fibered over tori, II},
 Proc. Japan Acad. Ser.A Math Sci. {\bf 59} (1983), 100-103.

 \bibitem{m}
 Y. Matsumoto,
 "Lefschetz fibrations of genus two - A topological approach",
 in {\em Topology and Teichm\"{u}ller Spaces (Katinkulta, Finland, 1995)},
 World Sci., River Edge, New Jersey. 1996, 123-148.

 \bibitem{mp}
 J. D. McCarthy and A. Papadopoulos,
 {\em Involutions in surface mapping class groups},
 Enseign. Math. (2) {\bf 33} (1987), 275-290.

 \bibitem{o}
 B. Ozbagci,
 {\em Signatures of Lefschetz fibrations},
 to appear in Pacific J. Math.

 \bibitem{os}
 B. Ozbagci and A. I. Stipsicz,
 {\em Noncomplex smooth $4$-manifolds with genus-$2$ Lefschetz fibrations},
 Proc. Amer. Math. Soc. {\bf 128} (2000), 3125-3128.

 \bibitem{s}
 A. I. Stipsicz,
 {\em On the number of vanishing cycles in Lefschetz fibrations},
 Math. Res. Lett. {\bf 6} (1999), 449-456.

 \bibitem{s2}
 A. I. Stipsicz,
 {\em Singular fibres in Lefschetz fibrations on manifolds with $b_2^+=1$},
 to appear in Topology Appl.

 \end{thebibliography}
 \end{document}